\documentclass[preprint,12pt]{elsarticle}
\usepackage{amsmath}
\usepackage{amssymb}
\usepackage{epsfig}
\usepackage{epstopdf}
\usepackage[latin5]{inputenc}
\usepackage{subfigure}
\usepackage{lineno,hyperref}
\usepackage[mathscr]{eucal}
\usepackage{color}
\usepackage{graphics}
\usepackage{csquotes}

\modulolinenumbers[5]

\usepackage{color}

\numberwithin{equation}{section}

\newtheorem{thm}{Theorem}

\numberwithin{equation}{section}
\numberwithin{thm}{section}
\numberwithin{lemma}{section}
\numberwithin{prop}{section}
\numberwithin{cor}{section}
\numberwithin{rmk}{section}
\numberwithin{defn}{section}

\definecolor{darkolivegreen}{rgb}{0.333333, 0.419608, 0.1843140}

 \DeclareMathOperator{\sn}{sn}

\DeclareMathOperator{\tn}{tn}

\newcommand{\dx}{\partial_x}

\newcommand{\dt}{\partial_t}
\newcommand{\du}{\partial_u}

\begin{document}

\title{On the Rosenau equation: Lie symmetries, periodic solutions and solitary wave dynamics}

\author[a]{Ali Demirci}
\ead{demircial@itu.edu.tr}

\author[b]{Yasin Hasano\u{g}lu}
\ead{y.hasanoglu2020@gtu.edu.tr}

\author[a]{Gulcin M. Muslu\corref{cor1}}
\ead{gulcin@itu.edu.tr}

\author[a]{Cihangir \"{O}zemir}
\ead{ozemir@itu.edu.tr}

\cortext[cor1]{Corresponding author}
\address[a]{Istanbul Technical University, Department of Mathematics, Maslak, Istanbul, Turkey}
\address[b]{Gebze Technical University, Department of Mathematics, Kocaeli, Turkey}

\begin{abstract}
In this paper,  we first consider the  Rosenau equation with the quadratic nonlinearity and identify its Lie symmetry algebra. We obtain reductions of the equation to ODEs,  and find periodic analytical solutions in terms of elliptic functions. Then,  considering a general power-type nonlinearity, we  prove the  non-existence of solitary waves \textcolor{black} {for some parameters} using Pohozaev type identities. The Fourier pseudo-spectral method is proposed for the
Rosenau equation with this  single power type nonlinearity. In order to investigate the solitary wave dynamics, we generate the solitary wave profile as an initial condition by using the Petviashvili's method. Then the evolution of the single solitary wave and overtaking collision of solitary waves are investigated by various numerical experiments.
\end{abstract}


\date{\today}

\maketitle

\section{Introduction}
The vibrations of a one dimensional anharmonic lattice associated with the birth of the soliton  are modeled in terms of the discrete lattices.
If the lattice is dense and weakly anharmonic, the well-known Korteweg-de Vries (KdV) equation is derived.
However, KdV equation  cannot model the wave to wave and wave to wall interactions for the dynamics of dense discrete systems.  To overcome this difficulty, the Rosenau equation
\begin{equation}
u_t+u_x+u_{xxxxt}+(u^2)_x=0 \label{rosenauquadratic}
\end{equation}
 is derived to describe the dynamics of dense discrete systems  considering higher order effects by Rosenau
\cite{rosenau1988dynamics}.

In this paper, we study the Rosenau equation with a single power-type nonlinearity
\begin{equation}\label{Rosenau}
u_t+u_x+u_{xxxxt} +(g(u))_x=0,
\end{equation}
where $g(u)=\displaystyle\frac{u^{p+1}}{p+1}$ with $p>0$. The Rosenau equation satisfies the conservation law \textcolor{black}{ \cite{Park1990}}:
\begin{equation}\label{energy}
\mathcal{E}=\int_\mathbb{R} \left[ u^2+({u_{xx}})^2 \right] ~dx.
\end{equation}

\noindent
Park \cite{Park1990} studied  the  global existence and uniqueness of the initial value problem
for the Rosenau equation of the form
\begin{eqnarray}
&& u_t+u_{xxxxt}=(f(u))_{x},    \hspace{20pt} x\in \mathbb{R},\,\, t>0 \label{ivp1}  \\
&& u(x,0)= \varphi (x),    \label{ivp2}
\end{eqnarray}
where $f(u)=\displaystyle\sum_{i=1}^{n}\frac{c_i}{p_i+1}\, u^{p_i+1}$ with $c_i\in \mathbb{R}$ and $p_i\in \mathbb{Z}^+$. It was proved that the Cauchy problem \eqref{ivp1}-\eqref{ivp2} admits a unique global solution for the initial data $\varphi \in H_0^4(\mathbb{R})$. Park also proved that the solution with small initial data decays like $t^{-1/5}$ in \cite{Park1992}. In the multidimensional case, the global existence and uniqueness result of the IBVP with Dirichlet boundary conditions for the equations
\begin{equation}
u_t+\sum_{i=1}^{n}c_i\frac{\partial^4}{\partial x_i^4}\,u_t =\nabla \cdot f(u)
\end{equation}
and
\begin{equation}\label{multiR}
u_t+\Delta^2 u_t=\nabla\cdot f(u)
\end{equation}
are given in \cite{Park1993}.

\par
From the numerical point of view,   the finite element Galerkin approximate solutions of  \eqref{rosenauquadratic}
and the  error estimates are studied  in  \cite{chungha}. Fully discrete schemes,  backward Euler, Crank-Nicolson and two step backward
methods, are proposed for  the initial-boundary value problem (IBVP) of the Rosenau equation
\begin{equation}\label{multiR}
u_t+\Delta^2 u_t=\nabla\cdot f(u), \quad  (x,t)\in \Omega\times(0,T_0)
\end{equation}
where $\Omega$ is a bounded set in $\mathbb{R}^d$, $d=1,2,3$  in \cite{Chunk2001}. The authors obtain optimal $L^2$ estimates using Galerkin approximations and derive \emph{a priori} estimates for discrete schemes.
When $d=1$ and $f(u)=-(u+u^2)$,  a discontinuous Galerkin method and stability analysis  are investigated in   \cite{Choo2008}.  The initial value problem \eqref{ivp1}-\eqref{ivp2} is solved in \cite{Danumjaya2019} via discontinuous Galerkin finite element methods. A  conservative unconditionally stable finite difference scheme is used for the equation for \eqref{rosenauquadratic} in \cite{Omrani2008}.  A second order splitting combined with orthogonal cubic spline collocation method is  employed in \cite{manickam}.  The authors of \cite{Erbay2021} consider a nonlocal nonlinear PDE, which reduce to Eq. \eqref{Rosenau} for a special case of the kernel function. \textcolor{black}{ They establish a numerical scheme based on truncated discrete convolution sums} applicable to  \emph{"geniunely nonlocal"} \cite{Erbay2021} case, where usual finite-difference schemes will not work.

In the  current work, we  analyze the Rosenau equation \eqref{Rosenau} through several different approaches. We start with the quadratic nonlinearity, namely with Eq. \eqref{rosenauquadratic}.  As the current literature does not contain any result regarding the Lie symmetry algebra of \eqref{rosenauquadratic}, we  identify the symmetry algebra in this case  and the equation is reduced to ordinary differential equations.
After this we aim at finding exact solutions of traveling wave type, which was achieved with obtaining several periodic solutions. Section 2 is devoted to this analysis.  The main outcome of Section 2 is the analytical solutions obtained in terms of the elliptic functions. Some of these solutions are smooth, and some have singularities. To the best of our knowledge, these solutions have been obtained the first time in the existing literature with this work.

The equations constructed by adding some terms to Rosenau equation  such as  Rosenau-KdV equation, Rosenau-Kawahara equation \cite{Zuo09} have exact solitary wave solutions. There have been many studies focusing on the solitary wave dynamics of these equations. To the best of our knowledge, there is no exact solitary wave solution for the Rosenau equation with single power type nonlinearity.
The solitary wave solution was first generated in \cite{erkip} by using Petviashvili's method, numerically.
The non-existence result of solitary wave solutions for some parameters is  given in Section $3$. We use Pohozaev type identities to show the non-existence of solitary wave.
The existence and stability of  the solitary wave solutions of the equation  \eqref{Rosenau} for $p<8$  is discussed in \cite{zeng}.   To the best of our knowledge, existence and stability of solitary waves for $p \geq 8$ is an open problem.  The present study addresses the existence of solitary wave  solutions for $p \geq 8$ numerically by using Petviashvili method. For the time evolution of the constructed solution we need an efficient numerical method. We therefore propose a numerical method combining a Fourier pseudo-spectral method for the space discretization and a fourth-order Runge-Kutta scheme for time discretization. Section 4 is devoted the fully-discrete Fourier pseudo-spectral scheme  and show how to formulate it for the Rosenau equation.  We also discuss the evolution of the single solitary wave solutions   and   interaction of two solitary waves for the Rosenau equation. As far as we know, the solitary wave dynamics of the   Rosenau equation with single power type nonlinearity
has never been investigated before in literature.  In Section 5, the numerical scheme is tested for accuracy and convergence rate. The solitary wave dynamics of the Rosenau equation by considering various problems, like propagation of single solitary wave,  collision of two solitary waves  are discussed.

\section{The Symmetry Algebra and Reductions of the Rosenau Equation}
In this Section we will determine the Lie symmetry algebra of \eqref{Rosenau} when $p=1$, i.e., in case of a quadratic nonlinearity. This invariance algebra is going to provide us several reductions of the equation to ordinary differential equations. One of these reductions, which determines the traveling solutions of \eqref{Rosenau}, will be solved in terms of elliptic functions.

When $p=1$, if we replace $u\rightarrow 2u$, \eqref{Rosenau} takes the form
\begin{equation}\label{main}
u_t+u_x+(u^2)_x+u_{xxxxt}=0.
\end{equation}
\textcolor{black}{Applying the well-known procedure for finding the Lie symmetries of a differential equation,  we obtain the following result.}
\begin{thm}
Lie symmetry algebra $\mathscr{L}$ of the Rosenau equation \eqref{main}
is 3-dimensional and is generated by the vector fields
\begin{equation}\label{algebra}
\mathscr{V}_1=\dt, \quad \mathscr{V}_2=t\dt-\left(\frac{1}{2}+u\right)\du, \quad \mathscr{V}_3 = \dx
\end{equation}
with the only nonzero commutation relation $[\mathscr{V}_1,\mathscr{V}_2]=\mathscr{V}_1$. The algebra has the direct sum structure $\mathscr{L}=\{\mathscr{V}_1,\mathscr{V}_2\}\oplus \mathscr{V}_3=\mathscr{A}_2\oplus \mathscr{A}_1$. The optimal system of one-dimensional subalgebras of $\mathscr{L}$ is given in \cite{patera1977subalgebras} as
\begin{equation}
\textcolor{black}{\mathscr{L}_1}=\{\mathscr{V}_1\}, \quad\textcolor{black}{ \mathscr{L}_2}=\{-\mathscr{V}_2 \cos \theta+\mathscr{V}_3 \sin \theta\}, \quad \textcolor{black}{ \mathscr{L}_3}=\{\mathscr{V}_1 + \epsilon \mathscr{V}_3\}
\end{equation}
where $\epsilon = \mp 1$ and $0 \leq \theta <\pi$.
\end{thm}

In the following we briefly study the reductions of \eqref{main} to ODEs making use of these one-dimensional subalgebras.

\emph{\textbf{(i)} \textcolor{black}{ Reduction through $\mathscr{L}_1$.}}  The solutions invariant under the transformations generated by  $\mathscr{V}_1$ have the form $u=u(x)$, hence we obtain the trivial constant solutions $u=C$.

\emph{\textbf{(ii)} \textcolor{black}{ Reduction through $\mathscr{L}_2$.}}   First suppose $\theta=\pi/2$. Then $-\mathscr{V}_2.0+\mathscr{V}_3 .1=\dx$, therefore one looks for solutions $u=u(t)$, which again gives the trivial solutions. Let $\theta \in [0,\pi)-\{\frac{\pi}{2}\}$.   In order to determine the solutions which are invariant under the group of transformations generated by \ $-\mathscr{V}_2 \cos \theta+\mathscr{V}_3 \sin \theta$, we solve the invariant surface condition
\begin{equation}
(-\mathscr{V}_2 \cos \theta+\mathscr{V}_3 \sin \theta)G(t,x,u)=0
\end{equation}
and find the two invariants $i_1=x-\tilde c\ln t$ and $i_2=t(\frac{1}{2}+u)$, where $\tilde c = -\tan \theta$. Considering $i_1$ as the new independent variable and $i_2$ as the new dependent variable, group-invariant solutions in this case will be searched according to $i_2=F(i_1)$ and hence
\begin{equation}
u(x,t)=-\frac{1}{2}+\frac{1}{t}F(x-\tilde c \ln t),
\end{equation}
of which replacement in \eqref{main} gives the reduction
\begin{equation}\label{red1}
\tilde c F^{(5)}+F^{(4)}+(\tilde c F-F^2)'+F=0.
\end{equation}
Observe that when $\tilde c=0$, the order of this reduction is one less:
\begin{equation}\label{red2}
u(x,t)=-\frac{1}{2}+\frac{1}{t}F(x), \qquad             F^{(4)}-(F^2)'+F=0.
\end{equation}

\emph{\textbf{(iii)} \textcolor{black}{Reduction through $\mathscr{L}_3$.}}  The reduction with \textcolor{black}{ the generator} $\mathscr{V}_1+\epsilon \mathscr{V}_3$ determines  the traveling wave solutions: When we plug $u(x,t)=F(x-\epsilon t)$
in \eqref{main} we get
\begin{equation}\label{red3}
\epsilon F^{(4)}-F^2+(\epsilon-1)F=K_1
\end{equation}
where $\epsilon = \mp 1$ and $K_1$ is a constant.

\textcolor{black}{In the following subsection, by proposing a suitable expansion for the solution, we are going to find some exact solutions for the reduction \eqref{red3} that yields traveling solutions. For the reductions  \eqref{red1} and \eqref{red2} similar treatments can be performed, which we kept beyond the scope of this article.}

\subsection{Elliptic type solutions}
We shall try to find some exact solutions of traveling wave type. We can work on \eqref{red3}, but we prefer to write down the reduction for the equivalent traveling wave ansatz: $u(x,t)=F(\xi)$, $\xi=kx-ct$. Plugging this in \eqref{main} and integrating twice  we obtain
\begin{equation}\label{eqF}
c k^4 \Big[ F'''F'-\frac{1}{2} (F'')^2\Big]-\frac{k}{3}F^3+\frac{c-k}{2}F^2+K_1F+K_0=0.
\end{equation}
For \eqref{eqF} we propose
\begin{subequations}\label{ansatz}
\begin{eqnarray}
     F(\xi)&=&a_0+a_1\varphi(\xi)+a_2 \varphi^2(\xi)+a_3 \varphi^3(\xi)+a_4 \varphi^4(\xi), \label{exF}\\
     \xi(x,t)&=&k x-ct,\\
     \Big(\frac{d \varphi}{d\xi}\Big)^2&=&c_0+c_1\varphi(\xi)+c_2 \varphi^2(\xi)+c_3 \varphi^3(\xi)+c_4 \varphi^4(\xi)=P(\varphi(\xi)). \label{eqfi}
\end{eqnarray}
\end{subequations}
\textcolor{black}{The expansion in \eqref{eqfi} is suggested so as to obtain elliptic function solutions and if possible the limiting cases trigonometric and hyperbolic ones. The form of the expansion \eqref{exF} is chosen so that a balancing between the terms is possible when inserted in the equation. In our case, it appears that this balance is only possible when the degree of the  expansion in \eqref{exF} is 4.  Upon this substitution, in the resulting expression we express all derivatives of $\varphi(\xi)$ in terms of $\varphi$ using \eqref{eqfi}. Afterwards, we look for the possibility that coefficients of $\varphi^j$, $j=0,1,2,...$ vanish.}

We assume
$c_1=c_3=0$ and $a_1=a_3=0$.
We find the following values for the remaining constants   $a_0$, $a_2$, $a_4$, $c_0$, $K_0$ and  $K_1$;
\begin{subequations}\label{par3}
	\begin{eqnarray}
	a_0&=&\frac{112cc_2^2k^4+c-k}{2k},\\
	a_2&=&560cc_2c_4k^3,\\
	a_4&=&840cc_4^2k^3,\\
	c_0&=&\frac{2c_2^2}{9c_4}, \label{eqc0}\\
	K_1&=&\frac{18ck-9k^2+c^2(41216c_2^4k^8-9)}{36 k}\\
	K_0&=&\frac{c}{8}+\frac{c^3}{24k^2}-\frac{c^2}{8k}-\frac{k}{24}-\frac{5152}{9}c^2c_2^4k^6(c-k)+\frac{2057216}{81}c^3c_2^6k^{10},
	\end{eqnarray}
\end{subequations}
where  $\varepsilon=\pm{1}$ and $ c_2 $ and $c_4$ are arbitrary. Let us note that numerical values of $K_0$ and $K_1$ will not be important in the solutions we shall determine below.

Now  we need to integrate \eqref{eqfi}, which takes the form
\begin{equation}\label{pol2}
\dot \varphi ^2=c_0+c_2 \varphi^2+c_4 \varphi^4=P(\varphi),
\end{equation}
and find $\varphi(\xi)$ and hence $u(x,t)$. Evaluation of the integral of  \eqref{pol2} depends on the factorization of the polynomial $P(\varphi)$.
Assume that $\varphi_1$, $\varphi_2$, $\varphi_3$ and $\varphi_4$  are roots of the equation $P(\varphi)=0$.
Let us call $\varphi^2 = z$. Then we need to solve
\begin{equation}\label{zsq}
c_4z^2+c_2z+c_0=0
\end{equation}
keeping in mind from \eqref{eqc0} that $c_0c_4=2c_2^2/9$.  The discriminant of this equation is $\Delta = c_2^2 / 9$.
There are two cases.

\vspace{.5cm}
\noindent \underline{\textbf{Case I.}  $ \Delta=0 $ } \quad  This implies that $c_2=c_0=a_2=0$. Eq. \eqref{eqfi} takes the form
\begin{equation}\label{int}
\frac{ d\varphi}{\sqrt{c_4\varphi^4}}=\epsilon d\xi,
\end{equation}  where  $\epsilon=\pm{1}$. Therefore we  find
\begin{equation}
\varphi(\xi)=\epsilon\frac{1}{\sqrt{c_4}(\xi-\xi_0)}.
\end{equation}The solution becomes
\begin{equation}\label{sech}
u(x,t)=\frac{c-k}{2k}+\frac{840c k^3}{(kx-ct-\xi_0)^4}.
\end{equation}
\vspace{.5cm}

\noindent \underline{\textbf{Case II.}  $ \Delta>0 $ } Then $\displaystyle z_1=-\frac{c_2}{3c_4}$ and  $\displaystyle z_2=-\frac{2c_2}{3c_4}$.
There are two possibilities for this case: If the signs of $c_2$ and $c_4$ are opposite then the equation  $P(\varphi)=0$ has four distinct real roots $\mp \sqrt{z_1}$, $\mp \sqrt{z_2}$, which we name as
\begin{equation}
\varphi_1=\Big(\frac{-2c_2}{3c_4}\Big)^{1/2}, \quad \varphi_2=\Big(\frac{-c_2}{3c_4}\Big)^{1/2}, \quad \varphi_3=-\Big(\frac{-c_2}{3c_4}\Big)^{1/2} \quad \varphi_4=-\Big(\frac{-2c_2}{3c_4}\Big)^{1/2}.
\end{equation}
If the signs of $c_2$ and $c_4$ are the same then the equation $\eqref{pol2}$ has four distinct complex roots $\mp \sqrt{|z_1|}\, i$, $\mp \sqrt{|z_2|}\,i.$

 Let $ c_4>0 $ and $ c_2<0 $ . In order that  \eqref{pol2} makes sense, the right hand side must be nonnegative. Thus we should consider the  intervals
\textbf{(a)} $\varphi>\varphi_1>\varphi_2>\varphi_3>\varphi_4$, \textbf{(b)} $\varphi_1>\varphi_2>\varphi>\varphi_3>\varphi_4$ and \textbf{(c)} $\varphi_1>\varphi_2>\varphi_3>\varphi_4>\varphi$ when integrating  \eqref{pol2}.

\vspace{.5cm}
\noindent \underline{\textbf{Case II.a}} \quad
Let us first write
\begin{equation}\label{int2}
\frac{ d\varphi}{\sqrt{c_4(\varphi-\varphi_1)(\varphi-\varphi_2)(\varphi-\varphi_3)(\varphi-\varphi_4)}}= \epsilon d\xi,
\end{equation}
where $\epsilon =\mp 1$. In the first hand, when  $\varphi>\varphi_1>\varphi_2>\varphi_3>\varphi_4$, using the results available in the handbook \cite{byrd2013handbook}, we  obtain
\begin{equation}\label{IIa}
\int_{\varphi_1}^{\varphi}
\frac{d\tau}{\sqrt{c_4(\tau-\varphi_1)(\tau-\varphi_2)(\tau-\varphi_3)(\tau-\varphi_4)}}
=\frac{1}{\sqrt{c_4}}\,g\sn^{-1}\Big(\sqrt{\frac{(\varphi_2-\varphi_4)(\varphi-\varphi_1)}{(\varphi_1-\varphi_4)(\varphi-\varphi_2)}},m\Big)
\end{equation}
for the integration of the left hand side of \eqref{int2}.
Using the substitution
\begin{equation}
\mathrm{sn}^2u = \frac{(\varphi_2-\varphi_4)(\tau-\varphi_1)}{(\varphi_1-\varphi_4)(\tau-\varphi_2)},
\end{equation}
where
\begin{equation}
g=\frac{2}{\sqrt{(\varphi_1-\varphi_3)(\varphi_2-\varphi_4)}}=2(\sqrt{2}-1)\sqrt{\frac{-3c_4}{c_2}}, \quad m^2=\frac{(\varphi_2-\varphi_3)(\varphi_1-\varphi_4)}{(\varphi_1-\varphi_3)(\varphi_2-\varphi_4)}
\end{equation}
for which we calculate that $m=2\sqrt{3\sqrt{2}-4}$, gives rise to  the elliptic function solution  $\varphi$ to \eqref{pol2}
\begin{equation}\label{varfi}
\varphi(\xi)=\frac{\varphi_1-\varphi_2R\mathrm{sn}^2\Big(\epsilon
\frac{\sqrt{c_4}}{g}(\xi-\xi_0),m\Big)}{1-R\mathrm{sn}^2\Big(\epsilon\frac{\sqrt{c_4}}{g}(\xi-\xi_0),m\Big)},
\end{equation} where $\displaystyle R=\frac{\varphi_1-\varphi_4}{\varphi_2-\varphi_4}=4-2\sqrt{2}$. Hence the solution to \eqref{main} can be written as follows
\begin{eqnarray}\label{sola}
u(x,t)=a_0&+&a_2\left[\frac{\varphi_1-\varphi_2R\mathrm{sn}^2\Big(\frac{\sqrt{c_4}}{g}(kx-ct-\xi_0),m\Big)}{1-R\mathrm{sn}^2\Big(\frac{\sqrt{c_4}}{g}(kx-ct-\xi_0),m\Big)}\right]^2 \nonumber\\
          &+&a_4\left[\frac{\varphi_1-\varphi_2R\mathrm{sn}^2\Big(\frac{\sqrt{c_4}}{g}(kx-ct-\xi_0),m\Big)}{1-R\mathrm{sn}^2\Big(\frac{\sqrt{c_4}}{g}(kx-ct-\xi_0),m\Big)}\right]^4.
\end{eqnarray}
\textcolor{black}{Since $\mathrm{sn}(-u)=-\mathrm{sn}(u)$, we removed $\epsilon$ in passing from \eqref{varfi} to \eqref{sola}. For the following  cases II.b, II.c, II.d and II.e there are similar calculations, therefore we skip the details of the calculations made by using the transformations available in \cite{byrd2013handbook}  and state only the results.}

\vspace{.5cm}
\noindent \underline{\textbf{Case II.b}}  \quad   When   $ c_4>0 $ and $ c_2<0 $, integrating for $\varphi_1>\varphi_2>\varphi>\varphi_3>\varphi_4$ we obtain
 the elliptic function solution  $\varphi$ to \eqref{pol2} as
\begin{equation}
\varphi(\xi)=\frac{\varphi_2-\varphi_1R\mathrm{sn}^2\Big(\epsilon\frac{\sqrt{c_4}}{g}(\xi-\xi_0),m\Big)}{1-R\mathrm{sn}^2\Big(\epsilon\frac{\sqrt{c_4}}{g}(\xi-\xi_0),m\Big)},
\end{equation}
where
\begin{equation}
R=2(\sqrt{2}-1), \qquad g=2(\sqrt{2}-1)\sqrt{\frac{-3c_4}{c_2}}, \qquad m=2\sqrt{3\sqrt{2}-4}.
\end{equation}
Hence the solution to \eqref{main} can be written as
\begin{eqnarray}\label{solb}
u(x,t)=a_0&+&a_2\left[\frac{\varphi_2-\varphi_1R\mathrm{sn}^2\Big(\frac{\sqrt{c_4}}{g}(kx-ct-\xi_0),m\Big)}{1-R\mathrm{sn}^2\Big(\frac{\sqrt{c_4}}{g}(kx-ct-\xi_0),m\Big)}\right]^2 \nonumber\\
          &+&a_4\left[\frac{\varphi_2-\varphi_1R\mathrm{sn}^2\Big(\frac{\sqrt{c_4}}{g}(kx-ct-\xi_0),m\Big)}{1-R\mathrm{sn}^2\Big(\frac{\sqrt{c_4}}{g}(kx-ct-\xi_0),m\Big)}\right]^4.
\end{eqnarray}

\vspace{.5cm}
\noindent \underline{\textbf{Case II.c}}  \quad   When  $ c_4>0 $ and $ c_2<0 $, integrating over $\varphi_1>\varphi_2>\varphi_3>\varphi_4>\varphi$ we obtain
 the solution to \eqref{main}  as
\begin{eqnarray}\label{solc}
u(x,t)=a_0&+&a_2\left[\frac{\varphi_4-\varphi_3R\mathrm{sn}^2\Big(\frac{\sqrt{c_4}}{g}(kx-ct-\xi_0),m\Big)}{1-R\mathrm{sn}^2\Big(\frac{\sqrt{c_4}}{g}(kx-ct-\xi_0),m\Big)}\right]^2 \nonumber\\
          &+&a_4\left[\frac{\varphi_4-\varphi_3R\mathrm{sn}^2\Big(\frac{\sqrt{c_4}}{g}(kx-ct-\xi_0),m\Big)}{1-R\mathrm{sn}^2\Big(\frac{\sqrt{c_4}}{g}(kx-ct-\xi_0),m\Big)}\right]^4
\end{eqnarray}
with
\begin{equation}
R=4-2\sqrt{2}, \qquad  g=2(\sqrt{2}-1)\sqrt{\frac{-3c_4}{c_2}}, \qquad m=2\sqrt{3\sqrt{2}-4}.
\end{equation}

Let $ c_4<0 $ and $ c_2>0 $. Since the right hand side of  \eqref{pol2}  must be nonnegative, we should consider the  intervals
\textbf{(d)} $\varphi_1>\varphi_2>\varphi_3>\varphi>\varphi_4$ and \textbf{(e)} $\varphi_1>\varphi>\varphi_2>\varphi_3>\varphi_4$ when integrating  \eqref{pol2}.

\vspace{.5cm}
\noindent \underline{\textbf{Case II.d}} \quad
In case $ c_4<0 $ and $ c_2>0 $, employing a suitable substitution on  $\varphi_1>\varphi_2>\varphi_3>\varphi>\varphi_4$
a  solution of \eqref{main} is found to be as follows
\begin{eqnarray}\label{sold}
u(x,t)=a_0&+&a_2\left[\frac{\varphi_4+\varphi_1R\mathrm{sn}^2\Big(\frac{\sqrt{-c_4}}{g}(kx-ct-\xi_0),m\Big)}{1+R\mathrm{sn}^2\Big(\frac{\sqrt{-c_4}}{g}(kx-ct-\xi_0),m\Big)}\right]^2 \nonumber \\
&+&a_4\left[\frac{\varphi_4+\varphi_1R\mathrm{sn}^2
\Big(\frac{\sqrt{-c_4}}{g}(kx-ct-\xi_0),m\Big)}{1+R\mathrm{sn}^2\Big(\frac{\sqrt{-c_4}}{g}(kx-ct-\xi_0),m\Big)}\right]^4,
\end{eqnarray}
where
\begin{equation}
 R=3-2\sqrt{2} \qquad  g=2(\sqrt{2}-1)\sqrt{\frac{-3c_4}{c_2}}, \qquad m=3-2\sqrt{2}.
 \end{equation}

\vspace{.5cm}
\noindent \underline{\textbf{Case II.e}} \quad  For $ c_4<0 $ and $ c_2>0 $, working with a substitution when $\varphi_1>\varphi>\varphi_2>\varphi_3>\varphi_4$  we obtain

\begin{eqnarray}\label{sole}
u(x,t)=a_0&+&a_2\left[\frac{\varphi_2+\varphi_3R\mathrm{sn}^2\Big(\frac{\sqrt{-c_4}}{g}(kx-ct-\xi_0),m\Big)}{1-R\mathrm{sn}^2\Big(\frac{\sqrt{-c_4}}{g}(kx-ct-\xi_0),m\Big)}\right]^2
\nonumber \\
          &+&a_4\left[\frac{\varphi_2+\varphi_3R\mathrm{sn}^2\Big(\frac{\sqrt{-c_4}}{g}
          (kx-ct-\xi_0),m\Big)}{1-R\mathrm{sn}^2\Big(\frac{\sqrt{-c_4}}{g}(kx-ct-\xi_0),m\Big)}\right]^4,
\end{eqnarray}
where
\begin{equation}
R=3-2\sqrt{2},\qquad  g=2(\sqrt{2}-1)\sqrt{\frac{-3c_4}{c_2}}, \qquad m=3-2\sqrt{2}.
\end{equation}

\vspace{.5cm}
\noindent \underline{\textbf{Case II.f}} \quad Now, we consider the case when the polynomial $P(\varphi)$ has four distinct complex zeros.
In order that  \eqref{pol2} makes sense, the right hand side must be nonnegative.\textcolor{black}{ Therefore, we should consider the case $ c_4>0 $ and $ c_2>0 $.
\eqref{pol2} takes the form
\begin{equation}
\dot\varphi^2=c_4(\varphi^2+\frac{c_2}{3c_4})(\varphi^2+\frac{2c_2}{3c_4}).
\end{equation}}
 Therefore we have
\begin{equation}
\int_{0}^{\varphi}
\frac{d\tau}{\sqrt{c_4(\tau^2+a^2)(\tau^2+b^2)}}
=\frac{1}{\sqrt{c_4}}\,g\tn^{-1}\Big(\frac{\varphi}{b},\frac{1}{\sqrt{2}}\Big)
\end{equation}
with the substitution
\begin{equation}
\mathrm{tn}^2u=\frac{\tau^2}{b^2}, \qquad  a^2=\frac{2c_2}{3c_4}, \qquad  b^2=\frac{c_2}{3c_4},  \qquad  g=\sqrt{\frac{3c_4}{2c_2}}.
\end{equation}
This gives rise to  the elliptic function solution  $\varphi$ to \eqref{pol2},
\begin{equation}
\varphi(\xi)=b\,\mathrm{tn}\Big(\epsilon \sqrt{\frac{2c_2}{3}}(\xi-\xi_0),\frac{1}{\sqrt{2}}\Big).
\end{equation}

Hence the solution to \eqref{main} can be written as follows
\begin{eqnarray}\label{solf}
u(x,t)=a_0&+&a_2\left[b\,\mathrm{tn}\Big(\epsilon \sqrt{\frac{2c_2}{3}}(kx-ct-\xi_0),\frac{1}{\sqrt{2}}\Big)   \right]^2  \nonumber\\
       &+&a_4\left[b\,\mathrm{tn}\Big(\epsilon \sqrt{\frac{2c_2}{3}}(kx-ct-\xi_0),\frac{1}{\sqrt{2}}\Big)   \right]^4.
\end{eqnarray}
\begin{figure}[ht!]
\centering
\subfigure[]{ 
\includegraphics[scale=0.45]{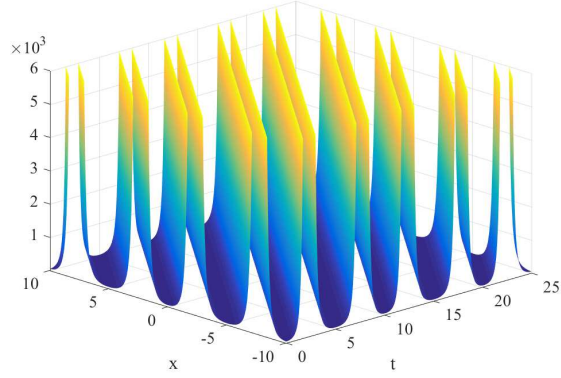}}
\subfigure[]{
\includegraphics[scale=0.45]{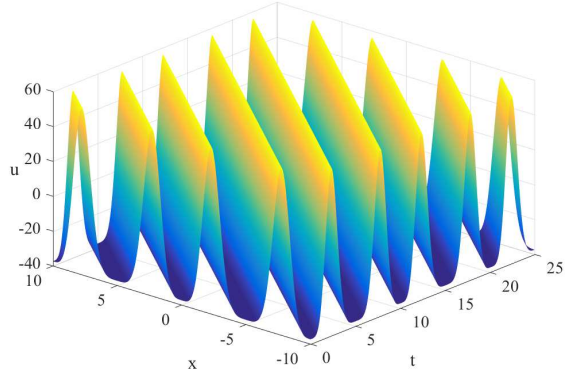}}
\caption{\small Plots of the periodic solution (a)$\eqref{sola}$,
(b) $\eqref{solb}$. For both of the solutions,  we choose the free parameters as $c=c_4=k=\epsilon=1$, $c_2=-1$ $\xi_0=0$. Note that there is a singular characteristic in the solution $\eqref{sola}$ whereas \eqref{solb} is smooth.}
\end{figure}

\textcolor{black}{We depict the solutions given in \eqref{sola} and \eqref{solb} for cases II.a and II.b, respectively. For both of the cases, we choose the arbitrary constants as}
\begin{equation}
c=c_4=k=1, \qquad c_2=-1, \qquad \xi_0=0.
\end{equation}
This determines the roots
\begin{equation}
\varphi_1=\sqrt{2/3}, \quad \varphi_2=1/\sqrt{3}, \quad \varphi_3=-1/\sqrt{3}, \quad \varphi_4=-\sqrt{2/3}.
\end{equation}
For both of the solutions, we have   $m=0.985171$,  $g=1.43488$, $a_0=56$, $a_2=-560$, $a_4=840$ and $c_0=2/9$.

For the solution \eqref{sola} of II.a, $R=4-2\sqrt{2}=1.17157$ for any choice of the arbitrary constants $c$, $c_2$, $c_4$ and $k$; therefore, \eqref{sola} has discontinuities. This is illustrated in Figure 1(a). On the contrary, for case II.b, $R=2(\sqrt{2}-1)=0.828427$ for any value of the arbitrary constants, therefore \eqref{solb} is smooth, which is well illustrated in Figure 1(b). Both solutions are periodic.

The other solutions \eqref{solc}, \eqref{sold}, \eqref{sole} and \eqref{solf} of the remaining cases give similar pictures to these two.

\setcounter{equation}{0}
\section{Solitary wave solutions}
In this section, we first establish the non-existence of solitary waves of the Rosenau equation   \eqref{Rosenau} \textcolor{black}{for some parameters.} To find the  localized solitary wave solutions of the Rosenau eq. \eqref{Rosenau}, we use the ansatz $u(x,t)=Q_c(\eta), ~~\eta=x-ct$
with $\displaystyle{\lim_{|\eta| \rightarrow \infty} Q_c(\eta)=0}$ which leads to  the ordinary differential equation
\begin{equation}
 -c Q_c^{\prime\prime\prime\prime\prime}  +(1-c) Q_c^{\prime}  + \frac{1}{p+1}(Q_c^{p+1})' =0.
  \label{solitary}
\end{equation}
Here $^\prime$ denotes the derivative with respect to $\eta$.
Integrating the equation \eqref{solitary} and, we have
\begin{equation}
c Q_c^{\prime\prime\prime\prime}+ (c-1) Q_c -\frac{Q_c^{p+1}}{p+1}=0.  \label{solitary2}
\end{equation}
The following theorem shows the non-existence of solitary waves \textcolor{black}{ for some parameters.}
\begin{thm}
The equation \eqref{solitary2} does not admit any nontrivial solution \mbox{$Q_c \in  H^{s}(\mathbb{R})$}  if one of the following conditions holds.
\begin{description}
\item[i.] $ c<0 $  and $p$ is even.
 \item[ii.]$ 0< c <1$,   for all $p>0$.
\end{description}
 \end{thm}
\textbf{Proof:} Let $Q_c$ be any nontrivial solution of the eq. \eqref{solitary2} in the class $H^{s}(\mathbb{R})$. Multiplying the eq. \eqref{solitary2} by $Q_c$, integrating on $\mathbb{R}$ and \textcolor{black}{ performing integration by parts twice for the first term,}  we get
\begin{equation}
   c\int_{\mathbb{{R}}} (Q_c^{\prime\prime})^2 dx + (c-1) \int_{\mathbb{{R}}} Q_c^2 dx = \frac{1}{p+1} \int_{\mathbb{{R}}} Q_c^{p+2} dx. \label{energy}
\end{equation}

\noindent
The term on the left side of this equation will be negative, a contradiction, when condition $(i)$ is satisfied.

\noindent
On the other hand, multiplying  the eq. \eqref{solitary2} by $x Q_c'$ and integrating over $\mathbb{R}$ yields the Pohozaev type identity
\begin{equation} \label{pohozaev}
  \frac{3c}{2}\int_{\mathbb{{R}}} (Q_c^{\prime\prime})^2 dx -\frac{c-1}{2}  \int_{\mathbb{{R}}} Q_c^2 dx
  = - \frac{1}{(p+1)(p+2)} \int_{\mathbb{{R}}} Q_c^{p+2} dx.
\end{equation}
Eliminating  $Q_c^{p+2}$ terms in the above equations gives
\begin{equation}
c\frac{3p+8}{2(p+2)}\int_{\mathbb{{R}}} (Q_c^{\prime\prime})^2 dx= (c-1) \frac{p}{2(p+2)}\int_{\mathbb{{R}}} (Q_c)^2 dx.
\end{equation}
The condition $0<c<1$ implies that the left hand side is non-negative and the right  hand side is negative. $\square$

\par
In \cite{zeng}, the author concerned with the class of  following equations
\begin{equation}
u_t+{\cal M}u_t +(f(u))_x=0,  \label{zengeq}
\end{equation}
where  $f$  is a real-valued function and ${\cal M}$ is a Fourier transform operator defined by
\begin{equation}
\widehat{{\cal M}u}(k)=m(k) \hat{u}(k)
\end{equation}
with $m(k)$ is an even and real valued function.
In Section $2$ of \cite{zeng},  existence and stability of solitary waves are proved
assuming  $ f(u)=u+\displaystyle\frac{u^{p+1}}{p+1} $ where $p>0$ is an integer and $p, m(k)$ satisfy the following conditions:
\begin{description}
	\item A1. there exist positive constants $A_1$ and $r>p/2$ \quad such that $m(k)\leq A_1|k|^r$ for $|k|\leq1$;
	\item A2. there exist positive constants $A_2$,$A_3$ and $s\geq1$ such that $A_2|k|^s\leq m(k)\leq A_3|k|^s$\quad for $|k|\geq1$;
	\item A3. $m(k)\geq0$ for all values of $k$;
	\item A4. $m(k)$ is four time differentiable for all non-zero values of $k$, and for each $j\in \{0,1,2,3,4\}$ there exist positive constants $B_1$ and $B_2$ such that
\end{description}
\begin{eqnarray*}
\Big\arrowvert\Big(\frac{d}{dk}\Big)^{j}\Big(\frac{m(k)-m(0)}{k}\Big)\Big\arrowvert&\leq& B_1|k|^{-j}\quad \text{for}\quad 0<|k|\leq1\\
\text{and}\qquad\qquad\qquad\qquad\qquad\qquad\qquad\qquad\qquad\qquad & & \nonumber\\ \Big\arrowvert\Big(\frac{d}{dk}\Big)^{j}\Big(\frac{\sqrt{m(k)}}{k^{{s}/{2}}}\Big)\Big\arrowvert&\leq& B_2|k|^{-j}\quad \text{for} \quad |k|\geq1.
\end{eqnarray*}

\noindent
 Choosing
${\cal M}=\partial_x^4$, the equation \eqref{zengeq} becomes the well-known Rosenau equation. The above theorem applies to the Rosenau equation in which  $m(k)=k^4$, when $p<8$. To the best of our knowledge, existence and stability of solitary waves for $p \geq 8$ is an open problem. In Section  \ref{sec:five}, we will answer this question numerically.

\setcounter{equation}{0}
\section{The numerical method}
We solve the Rosenau  equation by combining a Fourier pseudo-spectral method for the space component
and a fourth-order Runge Kutta scheme (RK4) for time.
If the spatial period $[a,b]$ is
normalized to $[0,2\pi]$ using the transformation
\mbox{$X=2\pi\displaystyle\frac{x-a}{b-a}$}, the equation \eqref{Rosenau} becomes
\begin{equation}\label{Rosenaut}
u_t+ \frac{2\pi}{b-a} ~ u_X+ (\frac{2\pi}{b-a})^4~  u_{XXXXt} + \frac{2\pi}{b-a} (\frac{u^{p+1}}{p+1})_X=0.
\end{equation}

\noindent
The interval $[0,2\pi]$ is divided into $N$ equal subintervals with grid spacing
$\Delta X=2\pi/N$, where the integer $N$ is even. The spatial grid points are given by
$X_{j}=2\pi j/N$,  $j=0,1,2,...,N$. The approximate solutions to
$u(X_{j},t)$ is denoted by $U_{j}(t)$.
The discrete Fourier transform of the sequence
$\{U_{j} \}$, i.e.
\begin{equation}\label{dft}
  \widetilde{U}_{k}={\cal F}_{k}[U_{j}]=
          \frac{1}{N}\sum_{j=0}^{N-1}U_{j}\exp(-ikX_{j}),
           ~~~~-\frac{N}{2} \le k \le \frac{N}{2}-1~
\end{equation}
gives the corresponding Fourier coefficients. Likewise, $\{U_{j} \}$ can
be recovered from the Fourier coefficients by the inversion formula
for the discrete Fourier transform (\ref{dft}), as follows:
\begin{equation}\label{invdft}
  U_{j}={\cal F}^{-1}_{j}[\widetilde{U}_{k}]=
          \sum_{k=-\frac{N}{2}}^{\frac{N}{2}-1}\widetilde{U}_{k}\exp(ikX_{j}),
          ~~~~j=0,1,2,...,N-1     ~.
\end{equation}
Here $\cal F$ denotes the discrete Fourier transform and
${\cal F}^{-1}$ its inverse. These transforms are  efficiently computed
using a fast Fourier transform (FFT) algorithm. In this study, we use  FFT
routines in Matlab (i.e. fft and ifft).

Applying the discrete Fourier transform to the equation \eqref{Rosenaut}, we obtain the first order ordinary differential equation
\begin{equation}
(\widetilde{U}_k)_t= -\frac{\displaystyle\frac{2\pi}{ b-a} ik}
{1+(\displaystyle\frac{2\pi k}{b-a})^4}\,[\,\widetilde{U}_k+\frac{(\widetilde{U^{p+1}})_k}{p+1}\,]. \label{rosenau-fourier}
\end{equation}
In order to handle the nonlinear term we use a pseudo-spectral approximation. We use the fourth-order Runge-Kutta method to solve the resulting ODE
\eqref{rosenau-fourier} in time.
Finally, we find the approximate solution by using the inverse Fourier transform \eqref{invdft}.
\setcounter{equation}{0}
\section{The numerical experiments} \label{sec:five}

In this section, we present some numerical experiments of the Fourier pseudo-spectral method for the Rosenau equation.  To the best of our knowledge, there is no exact solitary wave solution for the single power type nonlinearity $g(u)=\displaystyle\frac{u^{p+1}}{p+1},~ p>0$. In order to investigate the solitary wave dynamics, we first construct the solitary wave profile as an initial condition by using the Petviashvili's method. Then, by taking this initial condition the evolution of the single solitary wave and overtaking collision of solitary waves are investigated  by using Fourier pseudo-spectral method.
Since the exact solitary wave solution is unknown,  the "exact" solitary wave solution $u^{ex}$ is obtained numerically  with a very fine spatial step size  $N=1024$  and a very small time step $M=10000$   by using Fourier pseudo-spectral method.  In order to quantify the numerical results, the $L_\infty$-error norm is defined as
\begin{equation}
L_{\infty}\mbox{-error}=\max_i |~u^{ex}_i-U_i~|.
\end{equation}

\subsection{Accuracy test}
In order to test our scheme and to investigate dynamics of the solitary waves, we need an  initial condition. The initial condition is generated by using the Petviashvili's iteration  method  \cite{erkip, petviashvili, pelinovsky,yang}. The solitary wave solution   $u(x,t)=Q_c(x-ct)$  of the Rosenau equation satisfies the equation \eqref{solitary2}.
Applying   the Fourier transform to the equation \eqref{solitary2} yields
\begin{equation*}
( ck^{4}+ c-1 )  \widehat{Q_c}(k)=\frac{1}{p+1} \widehat{Q_c^{p+1}}(k).
\end{equation*}
The Petviashvili method  for the Rosenau eq. is given by
\begin{equation}
   \widehat{Q}_{n+1}(k)= (M_n)^{\nu} \frac{\widehat{Q_{n}^{p+1}}(k)}{(p+1)( ck^{4}+ c-1 ) } \label{scheme}
\end{equation}
with stabilizing factor
\begin{equation*}
  M_{n}=\frac{\int_{\mathbb{R}} [ck^4 +c-1] [\widehat{Q}_{n}(k)]^2 dk }{\frac{1}{p+1}\int_{\mathbb{R}}\displaystyle \widehat{Q}^{p+1}_{n}(k) \widehat{Q}_{n}(k)dk },
\end{equation*}
for some parameter $\nu$. Here $Q$ is used instead of $Q_c$ for simplicity. The Petviashvili's iteration method for Rosenau eq. was first introduced in \cite{erkip}. We refer to \cite{erkip} for detailed information.
The overall iterative process is  controlled by the  error,
\begin{equation}
  Error(n)=\|Q_n-Q_{n-1}\|,~~~~n=0,1,.... \nonumber
\end{equation}
 between two consecutive iterations defined with the  number of iterations,  the stabilization factor error
\begin{equation}
|1-M_n|, ~~~~n=0,1,.... \nonumber
\end{equation}
and the residual error
\begin{equation}
{RES(n)}= \|{\cal R} Q_n\|_\infty, ~~~~n=0,1,.... \nonumber
\end{equation}
where
\begin{equation}
{\cal R}Q= c Q^{\prime\prime\prime\prime} + (c-1) Q - \frac{1}{p+1} Q^{p+1}.
\end{equation}

\noindent
\textcolor{black}{The space interval is $-50 \le x \le 50$ and we choose the the number of spatial grid points $N=1024$.}
In the left panel of  Figure $2$, we present the solitary wave solution constructed by Petviashvili method with the speed $c=2$ for the quadratic nonlinearity. In the right panel of Figure $2$, we show the variation of    three different  errors
with the  number of iterations in semi-log scale.
\begin{figure}[!htbp]
\begin{minipage}[t]{0.45\linewidth}
\centering
\hspace*{-20pt}
\includegraphics[height=5.5cm,width=7.5cm]{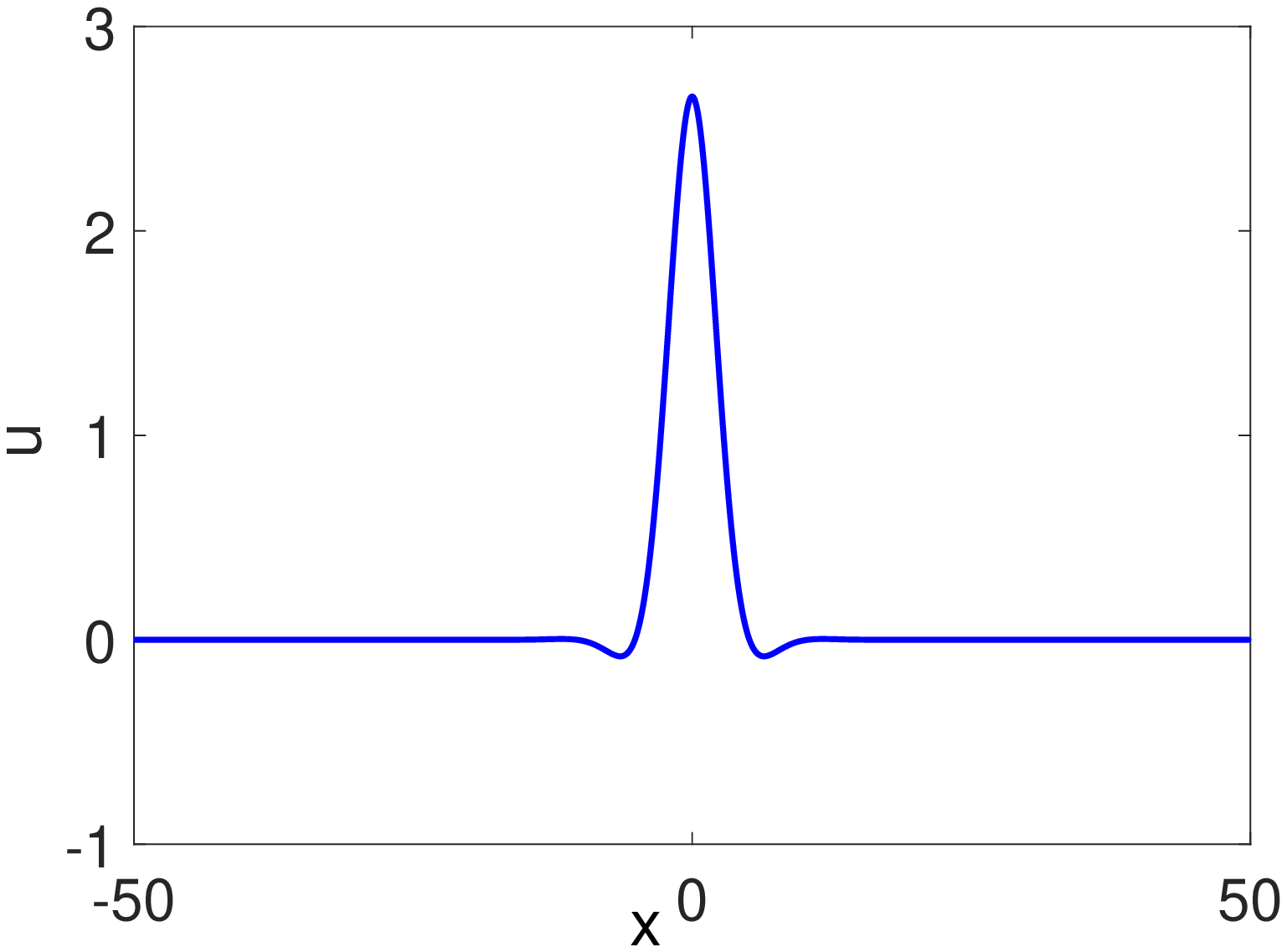}
\end{minipage}%
\hspace{3mm}
\begin{minipage}[t]{0.45\linewidth}
\centering
\includegraphics[height=5.5cm,width=7.5cm]{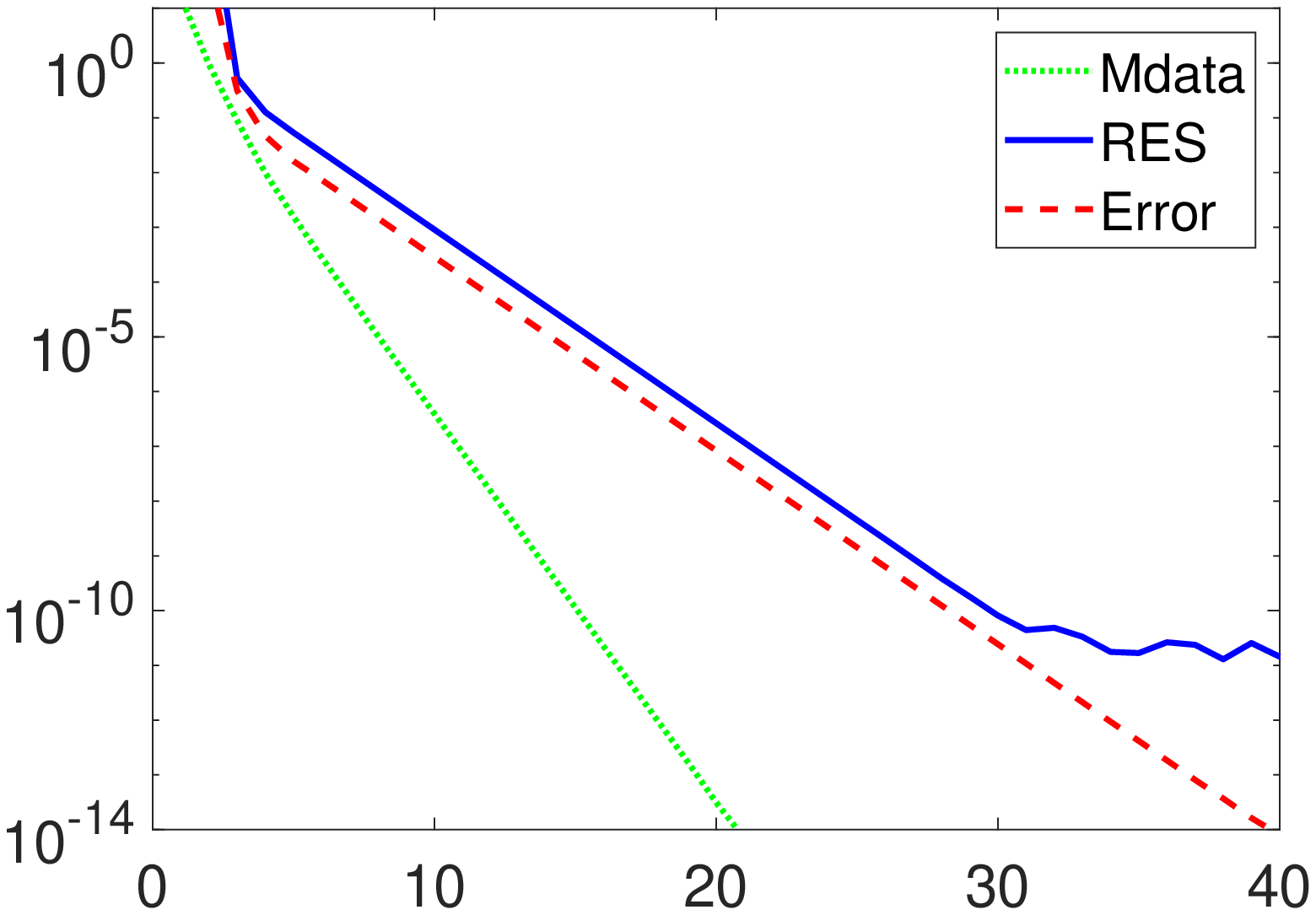}
\end{minipage}
\vspace{-5mm}
\caption{The solitary wave solution constructed by Petviashvili method with the speed $c=2$ for quadratic nonlinearity and  the variation of the $Error(n)$, $|1-M_n|$ and $RES$ with the number of iterations in semi-log scale.}\label{long}
\end{figure}
The solitary wave has an oscillatory structure.   In \cite{kawahara},  Kawahara  is concerned with the generalized KdV equation. When the coefficient of the fifth order derivative term is dominant over the third order derivative term, Kawahara observed the oscillatory solitary wave. The Rosenau equation does not have a third order term, the fifth order derivative term is dominant.
Therefore, this numerical  result is compatible with the result given in \cite{kawahara} and \textcolor{black}{ the result for the fifth order KdV equation given in  \cite{tranter}}.
The similar oscillatory wave for the wave speed $c=1.8$ is also observed in \cite{erkip}.

The solitary   wave profiles are constructed by using Petviashvili's  iteration  method for various nonlinearities.
\textcolor{black}{The problem is solved
$-50 \le x \le 50$ by taking $N=1024$. } In Figure 3, we only depict the solitary wave profiles for   $p=8,15$ and $p=30$
focusing \textcolor{black}{the interval $-20 \le x \le 20$.} The numerical experiments indicate that the solitary waves also exist for $p \ge 8$.  We observe that the amplitude of the solitary wave  decreases  with increasing nonlinear effects.
\begin{figure}[h!]
\centering
\includegraphics[height=7cm,width=7.5cm]{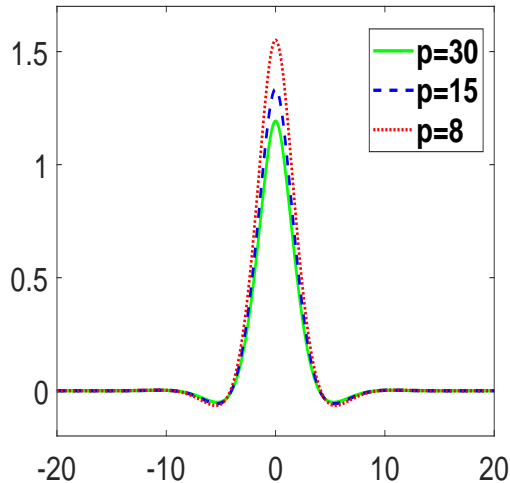}
\vspace{-5mm}
\caption{ The solitary wave solution constructed by Petviashvili method with the speed $c=2$ for $p=8,15$ and $p=30$.}\label{long}
\end{figure}

\begin{figure}[h!]
\centering
\includegraphics[height=6.5cm,width=8.5cm]{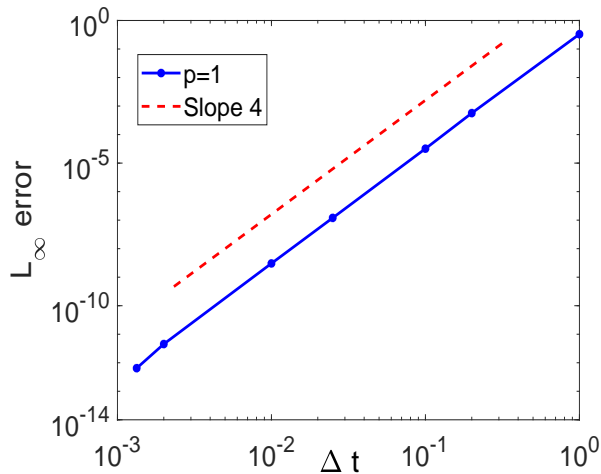}
\caption{The convergence rates in time calculated from
the $L_{\infty}$-errors. }\label{ordertime}
\end{figure}

\noindent
In order to check the temporal discretization errors, we fix the number of  spatial grid points $N=1024$ and solve the Rosenau eq. for different time step $\Delta t$. The convergence rates  calculated from the $L_{\infty}$-errors
at the terminating time $T=10$ are illustrated in Figure \ref{ordertime}. The computed convergence rates agree well
with the fact that Fourier pseudo-spectral method exhibits the  fourth-order convergence in time.
\noindent
In order to  test the spatial discretization errors,  we fix the the time step  such that the temporal error can be neglected, and solve the Rosenau eq. for different mesh size $\Delta x$.
In these experiments we take $M=10000$ to minimize  the temporal errors.  We  present the $L_{\infty}$-errors for the terminating time $T=10$
together with the observed rates of convergence in Table I. These results show that the numerical solution obtained using the
Fourier pseudo-spectral scheme converges rapidly to the accurate solution in space, which is an
indicative of exponential convergence.

\begin{table}
\caption{\label{tab:space}The convergence rates in space calculated from
the $L_{\infty}$-errors.}
\begin{center}
\begin{tabular}{|c|lc|} \hline\hline
$N$     & \hspace*{10pt}$L_{\infty}$-error \hspace*{5pt} & \hspace*{5pt} Order   \\
\hline
  32    &    0.891        & ~ -      \\
  64    &    1.277E-3     & ~9.4472  \\
  128   &    7.199E-10    & ~20.7581 \\
  256   &    1.998E-14    & ~15.1366 \\
\hline\hline
\end{tabular}
\end{center}
\end{table}
\noindent
\subsection{Single solitary wave}

\begin{figure}[ht]
 \begin{minipage}[t]{0.4\linewidth}
   \includegraphics[width=3.1in]{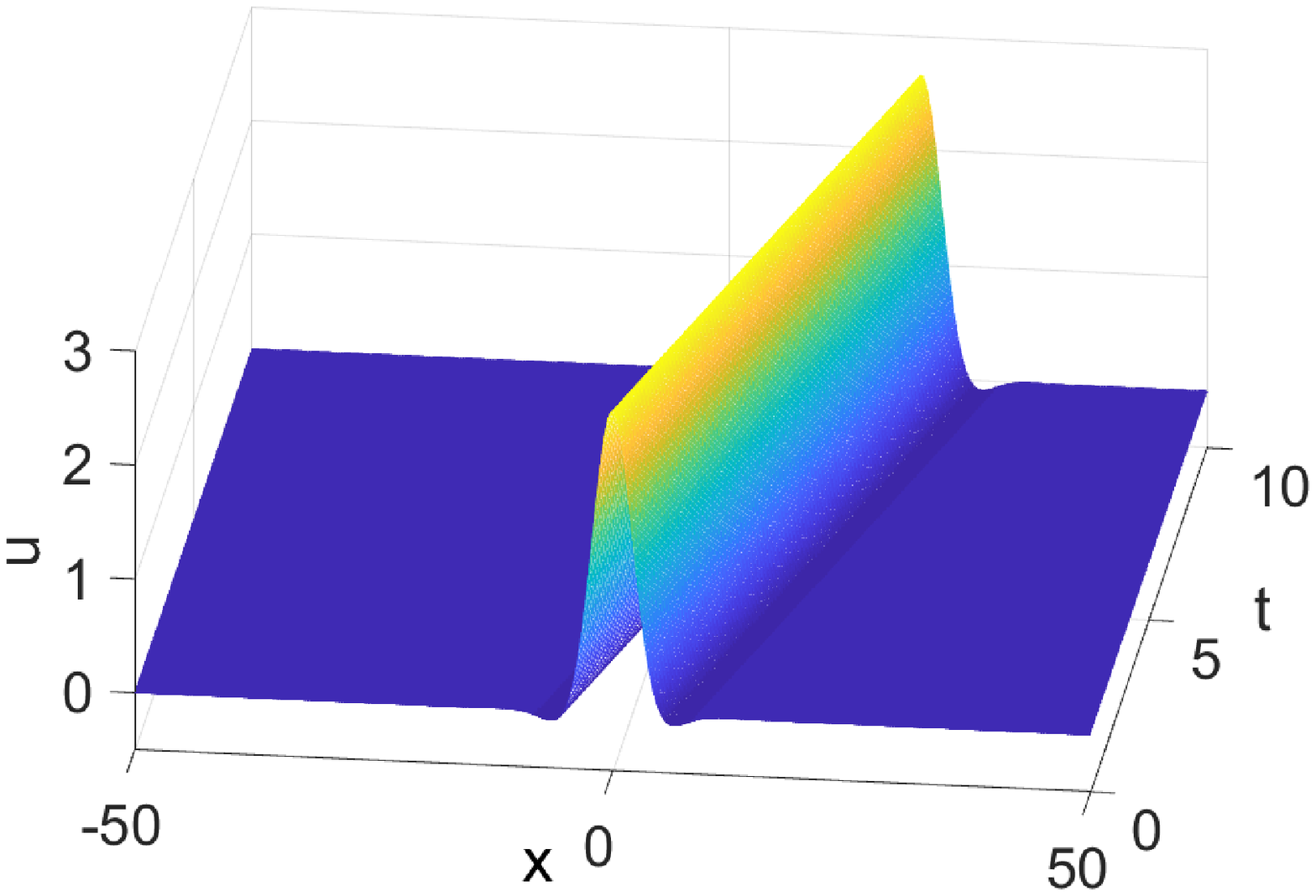}
 \end{minipage}
\hspace{50pt}
\begin{minipage}[t]{0.4\linewidth}
   \includegraphics[width=3in]{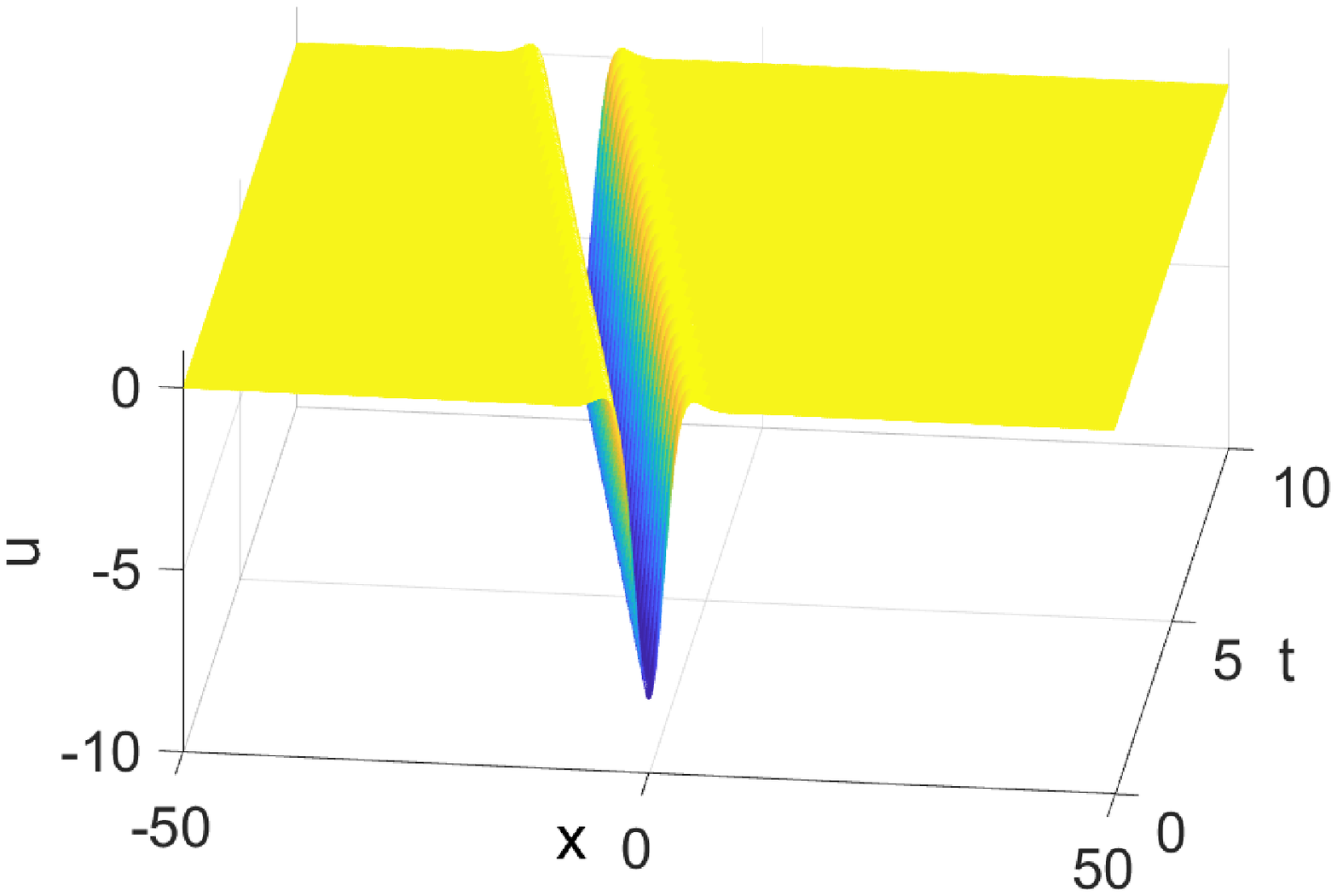}
 \end{minipage}
\caption{The evolution of the  solitary wave solution constructed by Petviashvili method with the speed $c=2$ (left panel) and $c=-2$ (right panel) by taking $N=1024$ and $M=10000$. \label{single}}
\end{figure}

In this subsection, we investigate time evolution of the numerically generated solitary waves by using a numerical scheme combining a Fourier pseudo-spectral method for space and a fourth order Runge-Kutta method for the time integration.
Computations are carried out with $N=1024$ and $M=10000$ on the interval $-50\le x\le 50$ \textcolor{black}{for times up to $T=10$.}   We present the evolution of the  solitary wave solution constructed by Petviashvili method with the speed $c=2$ and $c=-2$ for the quadratic nonlinearity in the left  panel and  in the right panel of  Figure \ref{single}, respectively. \textcolor{black}{ The maximum change in  the energy  $|{\cal E} (t) - {\cal E} (0)|$ is approximately $7.10 \times 10^{-14}$   and  $6.25 \times 10^{-13}$  during the entire time interval $0\le t \le 10$ for $c=2$ and $c=-2$, respectively.} This behavior provides a valuable check on the
numerical results.


\subsection{Interaction of two solitary waves}

\begin{figure}[!htbp]
 \begin{minipage}[t]{0.5\linewidth}
   \includegraphics[width=3.1in]{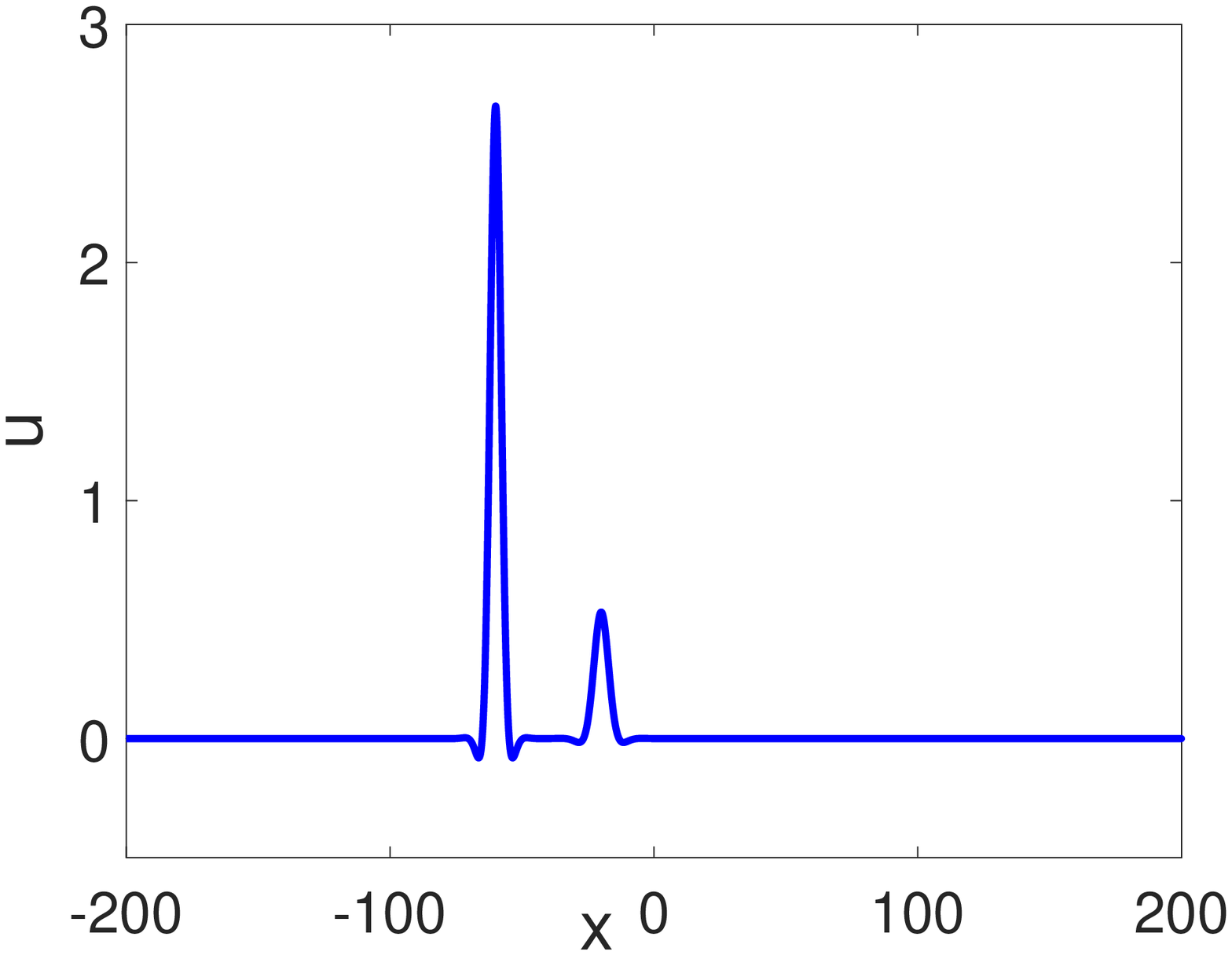}
 \end{minipage}
\hspace{20pt}
\begin{minipage}[t]{0.4\linewidth}
   \includegraphics[width=3in,height=2.5in]{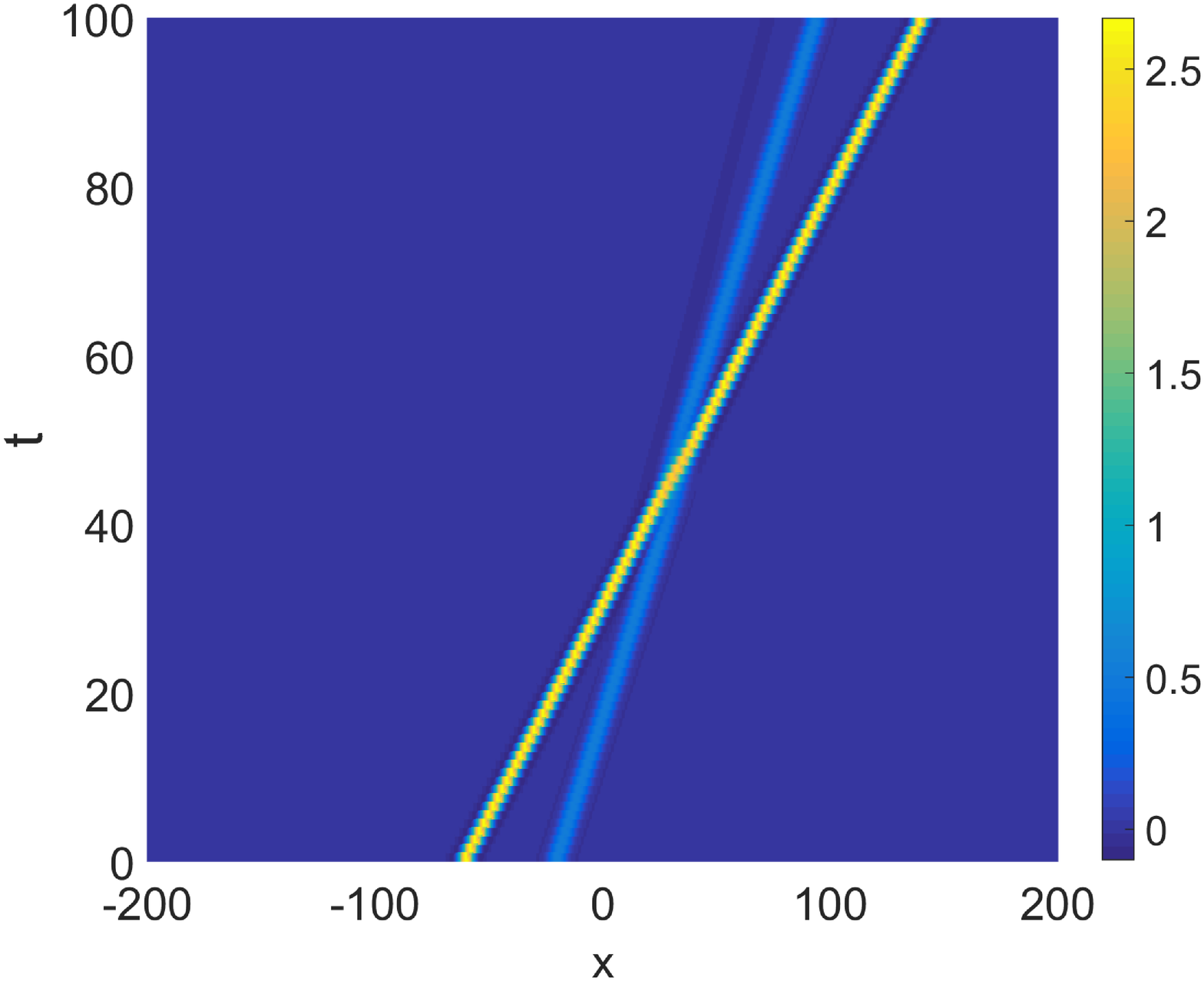}
 \end{minipage}
 \caption{The generated two solitary waves with speeds $c_1=2$ and $c_2=1.2$ by using Petviashvili's method (left panel) and surface plot of the interaction of two solitary waves (right panel) ($N=2^{14}$ and $M=10000$) \label{collision1} }
\end{figure}

In this subsection, we present some results to investigate the interaction of two solitary waves for quadratic nonlinearity.
The experiments in this subsection are performed  over the interval  $[-200, 200]$ with $N=2^{14}$ and $M=10000$ \textcolor{black}{for times up to $T=100$.}
The left panel of Figure \ref{collision1} shows the  generated two solitary waves   which are initially located at the positions $x_1=-60$ and $x_2=-20$ moving in the same direction  with speeds $c_1=2$ and $c_2=1.2$ by using Petviashvili's method, respectively.

In the right panel of
Figure \ref{collision1},  we illustrate the interaction of two solitary waves. It can be observed that the faster wave overtakes the slower one at around $t=50$  and leaves it behind as time evolves. We observe that the amplitude of the solution at the collision time   is smaller than the summation of the two initial amplitudes.

\begin{figure}[h!bt]
 \begin{minipage}[t]{0.5\linewidth}
   \includegraphics[width=3.1in,height=2.8in]{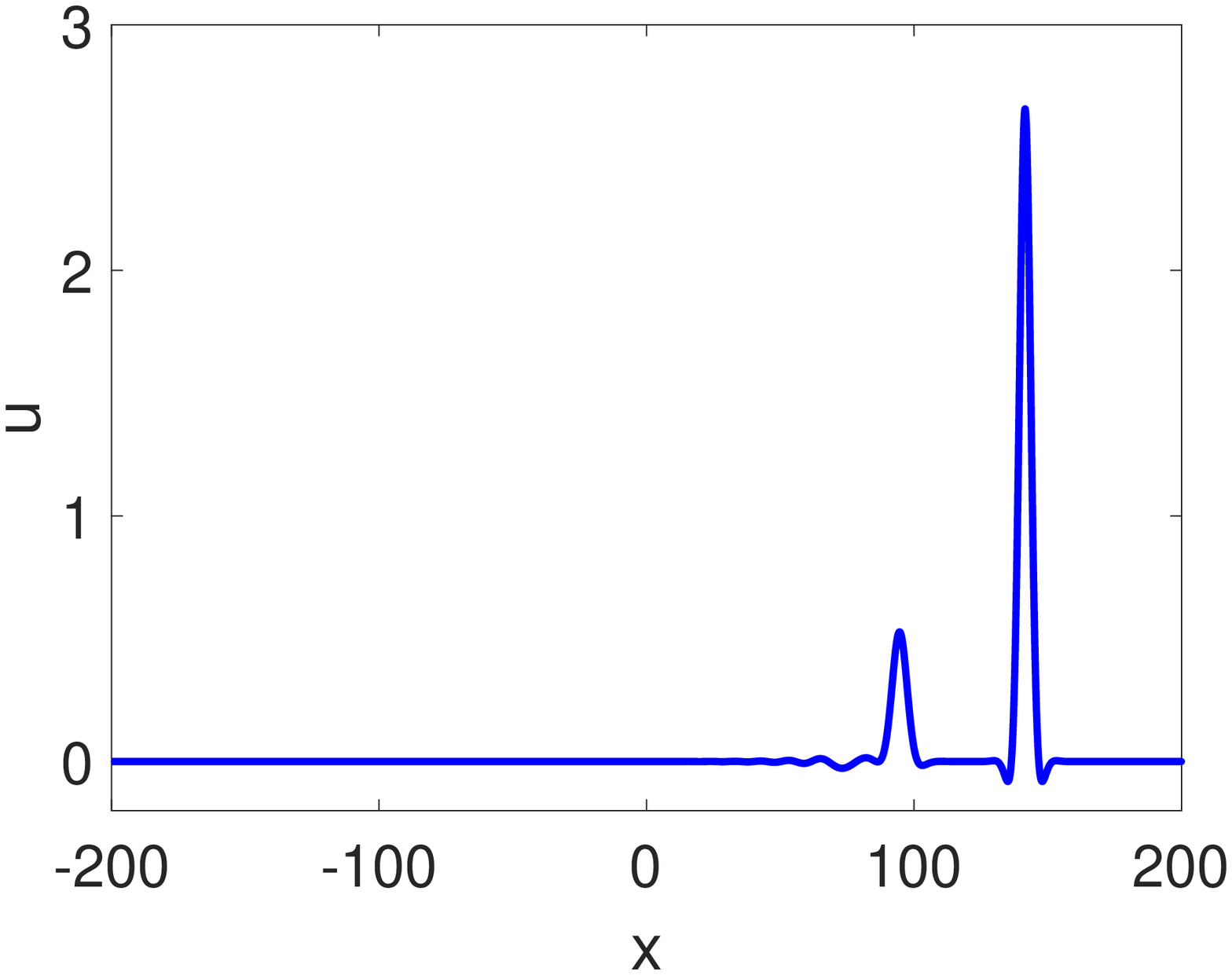}
 \end{minipage}
\hspace{20pt}
\begin{minipage}[t]{0.4\linewidth}
   \includegraphics[width=3in,height=2.8in]{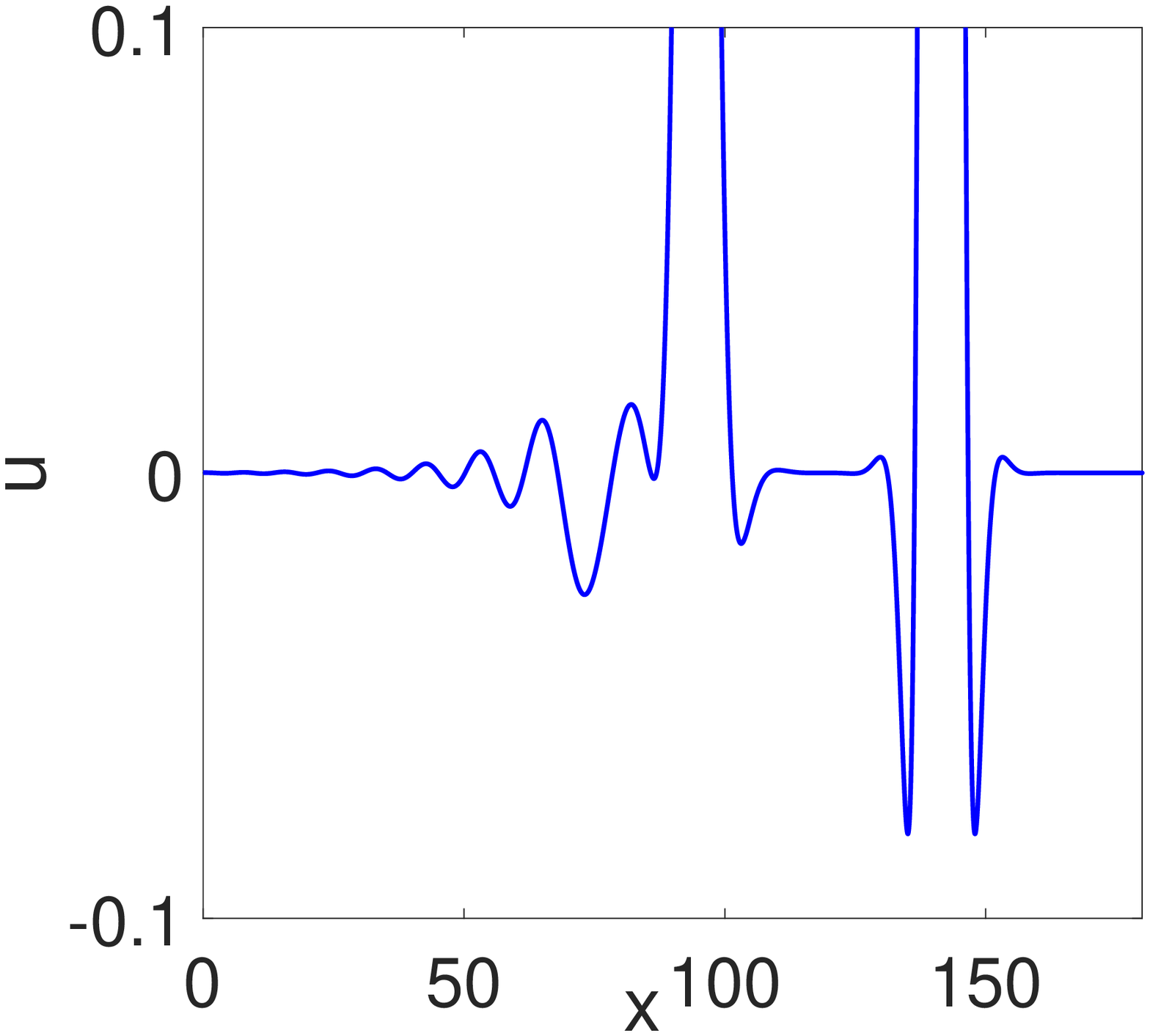}
 \end{minipage}
 \caption{The solution profile at the final time $T=100$ (left panel)  and the close-up look of the oscillatory tail (right panel)   \label{tail}}
\end{figure}
\begin{figure}[h!]
\centering
\includegraphics[height=6.5cm,width=8.5cm]{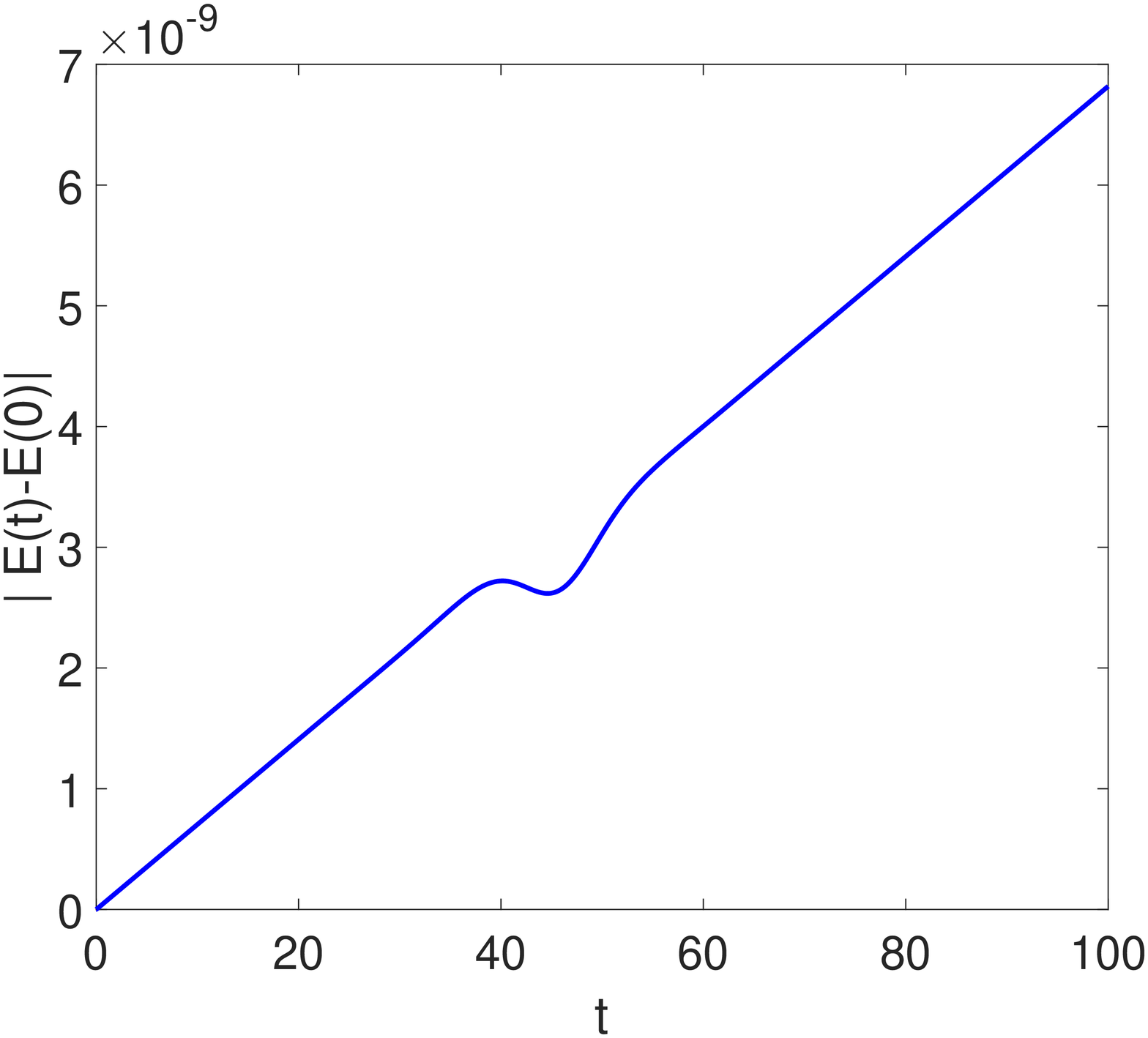}
\caption{The variation of the change  in  the conserved quantity  $\cal E$. }
\label{energy1}
\end{figure}

\textcolor{black}{The left panel of Figure \ref{tail} shows the solution at the final time  $T=100$. We observe the oscillatory wave trailing behind  the smaller wave after the interaction. The right panel of Figure \ref{tail} gives a closer look at the oscillatory tail.
To check the accuracy of the numerical solution in the dispersive tail region as in \cite{bona},  we obtained the numerical solutions taking $M=10000$ and $M=20000$. We compare the  dispersive tail of the numerical solutions correspond to   $M=10000$  and $M=20000$  with the dispersive tail of the numerical solution obtained by  $M=40000$. The $L_\infty$-error norm is $2.2946 \times 10^{-10}$ and $1.3662 \times 10^{-11}$, respectively. The ratio of the errors is approximately $16.8$ which is compatible with the fourth order convergence of the Runge-Kutta time stepping scheme.
Since the presence of dispersive tail, the collision between two solitary waves for the Rosenau equation is inelastic which indicates that the
Rosenau equaion is not integrable. We also present the
evolution of the  change  in  the conserved quantity  $\cal E$ (energy)   in  Figure \ref{energy1}.
It can be clearly observed that the Fourier pseudo-spectral method preserves the energy $\cal E$  very well.}
\newpage
\setcounter{equation}{0}
\section{Conclusion} \label{sec:six}
In this work we have considered  the Rosenau equation with a general power-type nonlinearity. We first took a group-theoretical point of view and presented the Lie algebra of the equation in the widely considered quadratic case. We found reductions of the equation to ODEs. Among these, we focused on the reduction that produces traveling type solutions. We successfully  obtained some analytical solutions in terms of elliptic functions.

After this achievement, taking another approach, we have proved the non-existence of the solitary waves for some parameters by using Pohozaev identities. We have used the Petviashvili's method to construct the solitary wave profile numerically.  Existence of solitary waves  for the Rosenau equation with the nonlinearity  $p \geq 8$ is  an open problem. The numerical experiments indicate that the solitary wave solution also exists for $p \geq 8$.  We have proposed a numerical method combining a Fourier pseudo-spectral method for the space discretization and a fourth-order Runge-Kutta scheme for time discretization for the  Rosenau equation.  The method  converges  the fourth order in time and spectrally in space.
To the best of our knowledge,  solitary wave dynamics for the Rosenau equation with single power type nonlinearity is investigated the first time with this work in the literature. Some interesting numerical experiments such as evolution of single solitary wave and interaction of solitary waves are performed by various numerical experiments.

\vspace*{20pt}
\noindent \textbf{Acknowledgements}\\
The authors would like to express  sincere gratitude to the reviewers for their constructive suggestions which helped to improve the quality of this paper.


\bibliographystyle{elsarticle-num}

\begin{thebibliography}{9}

\bibitem{rosenau1988dynamics}
P. Rosenau, Dynamics of dense discrete systems: High order effects, Prog. Theor. Phys. 79(5) (1988) 1028--1042.

\bibitem{Park1990}
M.I. Park, On the Rosenau equation, Mat. Aplic. Comp. 9(2) (1990) 145--152.

\bibitem{Park1992}
M.I. Park, Pointwise decay estimates of solutions of the generalized Rosenau equation, J. Korean Math. Soc. 29(2) (1992) 261--280.

\bibitem{Park1993}
M.I. Park, On the Rosenau Equation in multidimensional space, Nonlinear Anal.-Theo. 21(1) (1993) 77--85.
\bibitem{chungha}
 S. K. Chung ,  S. N. Ha. Finite element Galerkin solutions for the Rosenau equation, Appl. Anal. 54.1-2 (1994)  39-56.

\bibitem{Chunk2001}
S.K. Chunk, A. K. Pani, Numerical methods for the Rosenau equation, Appl. Anal. 77(3-4) (2001) 351--369.

\bibitem{Choo2008}
S.M. Choo, S.K. Chong, K.I. Kim, A discontinuous Galerkin method for the Rosenau equation, Appl. Numer. Math. 58(6) (2008) 783--799.

\bibitem{Danumjaya2019}
P. Danumjaya, K. Balaje, Discontinuous Galerkin Finite Element Methods for 1D Rosenau Equation, arXiv preprint arXiv:1911.12795 (2019).

\bibitem{Omrani2008}
K. Omrani, F. Abidi, T. Achouri, N. Khiari, A new conservative finite difference scheme for the Rosenau equation, Appl. Math. Comput. 201(1-2) (2008) 35--43.

\bibitem{manickam}
S. A. Manickam,  A. K. Pani,  S. K. Chung. A second order splitting combined with orthogonal cubic spline collocation method for the Rosenau equation, Numer. Methods Partial Differen. Equat. 14.6 (1998) 695-716.

\bibitem{Erbay2021}
H. A. Erbay, S Erbay, A. Erkip, A semi-discrete numerical method for convolution-type unidirectional wave equations, J. Comput. Appl. Math. 387 (2021) 112496.

\bibitem{Zuo09}
J.M. Zuo, Solitons and periodic solutions for the Rosenau-KdV and Rosenau-Kawahara equations, Appl. Math. Comput. 215(2) (2009) 835-840.

\bibitem{Zhou2016}
D. Zhou, C. Mu, Homogeneous initial-boundary value problem of the Rosenau equation posed on a finite interval, Appl. Math. Lett. 57 (2016) 7--12.



\bibitem{Safdari-Vaighani2018}
A. Safdari-Vaighani, E. Larsson, A. Heryudono, Radial Basis Function Methods for the Rosenau Equation and Other Higher Order PDEs, J. Sci. Comput. 75(3) (2018) 1555--1580.

\bibitem{erkip}
H. A. Erbay, S Erbay, A. Erkip, Numerical computation of solitary wave solutions of the Rosenau equation, Wave Motion 98 (2020) 102618.



\bibitem{patera1977subalgebras}
J. Patera, P. Winternitz, Subalgebras of real three-and four-dimensional Lie algebras, J. Math. Phys. 18(7) (1977) 1449--1455.


\bibitem{byrd2013handbook}
P.F.  Byrd, M.D. Friedman, Handbook of elliptic integrals for engineers and physicists, Springer, 2013.

\bibitem{zeng}
L. Zeng, Existence and stability of solitary wave solutions of equations of Benjamin Bona Mahony type, J. Differ. Equations 188 (1) (2003) 1--32.





\bibitem{petviashvili}
V. I. Petviashvili, Equation of an extraordinary soliton, Fiz. Plazmy. Phys. 2 (1976) 469--472.

\bibitem{pelinovsky}
D. E. Pelinovski, Y. A. Stepanyants, Convergence of {P}etviashvili's iteration method for numerical approximation of stationary solutions of nonlinear wave equations, SIAM J. Numer. Anal. 42 (2004) 1110-1127.

\bibitem{yang}
 J. Yang, Nonlinear waves in integrable and nonintegrable systems, Society for Industrial and Applied Mathematics, 2010.

\bibitem{kawahara}
T. Kawahara, Oscillatory solitary waves in dispersive media, J. Phys. Soc. Jpn. 33(1) (1972) 260--264.

\bibitem{tranter}
\textcolor{black}{K. R. Khusnutdinova, Y. A. Stepanyants,  M. R. Tranter. Soliton solutions to the fifth-order {K}orteweg-{d}e Vries equation and their applications to surface and internal water waves, Phys.  Fluids 30(2) (2018) 022104.}

\bibitem{bona}
\textcolor{black}{J. L. Bona,  H. Kalisch, Models for internal waves in deep water, Discrete Contin. Dyn. Syst. 6(1) (2000) 1--20.}

\end{thebibliography}
\biboptions{compress}

\end{document}